\documentclass[fleqn,10pt]{SelfArx} 

\usepackage[english]{babel} 

\usepackage{lipsum} 

\definecolor{color1}{RGB}{0,0,90} 
\definecolor{color2}{RGB}{0,20,20} 

\usepackage{hyperref} 
\hypersetup{hidelinks,colorlinks,breaklinks=true,urlcolor=color2,citecolor=color1,linkcolor=color1,bookmarksopen=false,pdftitle={Title},pdfauthor={Author}}

\usepackage{graphicx,epsfig,graphics,epsf}
\usepackage{float,amsmath,amsfonts,amssymb,amsthm}
\usepackage{algorithm}
\usepackage{algpseudocode}
\usepackage{subcaption}

\newtheorem{thm}{Theorem}
\newtheorem{lem}{Lemma}
\newtheorem{definition}{Definition}

\usepackage{bm}

 \newcommand{\bx}{\mathbf{x}}

\newcommand{\bw}{\mathbf{w}}

\newcommand{\bz}{\mathbf{z}}

\newcommand{\bv}{\mathbf{v}}
\newcommand{\bq}{\mathbf{q}}
\newcommand{\br}{\mathbf{r}}

\JournalInfo{{\includegraphics[width=0.25\textwidth,height=10.0mm]{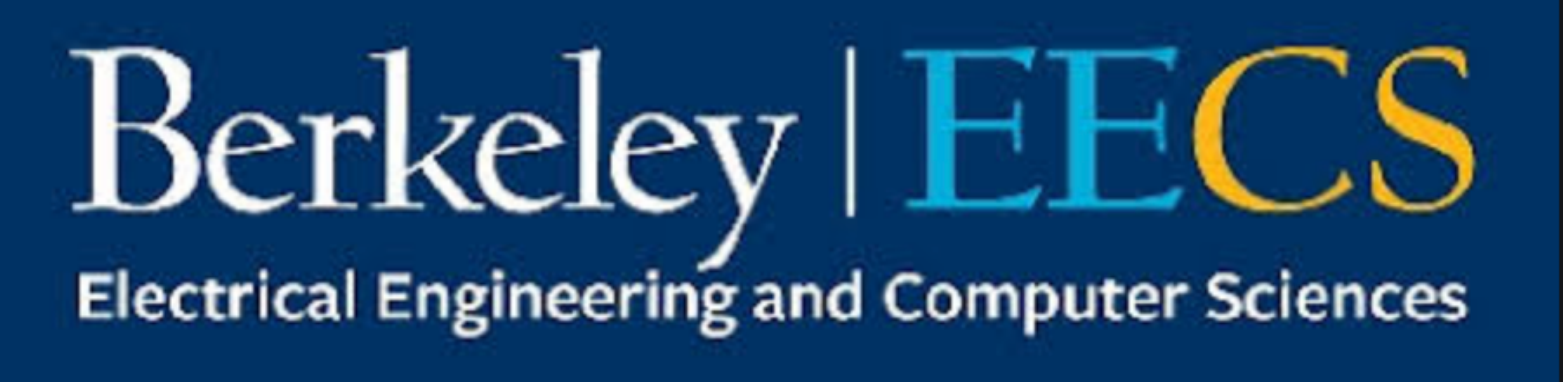}}\\Draft, 8/28/2020} 
\Archive{\fancyhead[R]{\includegraphics[width=0.15\textwidth,height=6.0mm]{Logo.pdf}}} 

\PaperTitle{Convergence of a Distributed Kiefer-Wolfowitz Algorithm} 

\Authors{Jean Walrand, Dept. of EECS, University of California, Berkeley} 
\affiliation{walrand@berkeley.edu} 

\newcommand{\keywordname}{Keywords} 

\onecolumn

\Abstract{This paper proposes a proof of the convergence of a distributed and asynchronous version of the Kiefer-Wolfowitz algorithm.}

\onecolumn

\begin{document}
\onecolumn

\flushbottom 

\maketitle 


\thispagestyle{empty} 


\section{Introduction}

The goal is to maximize a concave function of $K > 1$ variables.  There are $K$ agents and each agent observes the values of the function,
corrupted by observation noise,
and adjusts his own variable without knowing the values of the other variables.  The agents do not communicate their variables.
This formulation is motivated by many applications where the agents do not know each other or are not be able to 
communicate directly with one another.
Moreover, the agents are not synchronized, so that they 
update their variable either at the same or different times.

Each agent {\em experiments} by perturbing his variable by a zero-mean change in order to estimate the partial derivative of the 
function with respect to that variable.
He then {\em updates} his variable in proportion to the estimate of the partial derivative.  

This algorithm is an extension of
 \cite{KW52} and \cite{S92}. In \cite{KW52}, the authors introduce a gradient descent algorithm where the
gradient is estimated by observing the  function at perturbed values of its variable and they prove the convergence
of the algorithm to the minimum of the function.  \cite{S92} proposes a variation of the algorithm in the multivariate
case where the partial derivatives with respect to the different variables are estimated by simultaneously perturbing each variable by an independent
and zero-mean amount, instead of perturbing the variables one at a time.  
The author proves the convergence to the minimum of the function under the assumption that the variables
return infinitely often to a compact set.  In this paper, we extend the algorithm to the case where the different variables get updated
asynchronously.   Also, the proof does not require assuming returns to a compact set.   

An agent corrupts the estimate of the partial derivative of
another agent either when he experiments or updates his variable while the other agent calculates his estimate. 
Technically, the difficult aspect of our version is that
the corruption of the estimate by the updates of other agents is not zero-mean, in contrast with the corruption by their experiments
which is zero-mean.  Proving the convergence of this asynchronous version requires careful bounds on the size of the corruptions.
This is the technical contribution of the paper.

Some papers propose mechanisms where agents exchange the value of their variables, possibly with some delays, and they may know the function they want to maximize
(e.g., \cite{K19}, \cite{N09}, \cite{R19}, \cite{S20}).  The key contribution of this paper is to show that such communication is not necessary for convergence.
Also, the agents observe the values of the function with some observation noise but need not know its functional form. That is, the agents can observe the 
effect of their choice of value for their variable, but they could not calculate it.  The algorithm is similar  in spirit to t\^{a}tonnement (groping) in economics (see \cite{W1874}).

\section{Algorithm and Result}

Let $f: \Re^K \to \Re$ be a concave function, strictly concave in a neighborhood of its maximizer $\bx^*$.  
Assume that the function is globally Lipschitz with constant $L$.  Assume also that the second and third derivatives
of $f(\cdot)$ are bounded. Let $\tau \geq 2$ be an integer and $p_k \in \{0, \ldots, \tau - 1\}$ for $k = 1, \ldots, K$.
Let also $T_k(n) = n\tau + p_k - \tau 1\{p_k = \tau - 1\}$ for $k = 1, \ldots, K$ and $n \geq 0$.  Note that
$T_k(n) + 1 \in \{n \tau, \ldots, n\tau + \tau - 1\}$.
For $k = 1, \ldots, K$, agent $k$ {\em experiments} at times $\{T_k(n) + 1, n \geq 0\}$
and {\em updates} at times $\{T_k(n) + 2, n \geq 0\}$. Thus, the agents experiment and update every $\tau$
steps and they may be out of phase with one another. 
The case  of a single agent (i.e., $K = 1$) is the same as in \cite{KW52} while
that of simultaneous updates (i.e., $p_1 = \cdots = p_K$) corresponds to \cite{S92}.

The experiments and updates are defined as follows.
For $k = 1, \ldots, K$ and $n \geq 0$, let $x_k(n)$ be the value of the variable of agent $k$ at step $n$. 
Let also $\bx(n)$ be the vector with components $x_k(n)$, for  $n \geq 0$.

The algorithm is as follows.  For $k = 1, \ldots, K$ and $n \geq 0$, one has, for $m = T_k(n)$,
\begin{eqnarray}
&& x_k(m+1) = x_k(m) + a_k(n) \epsilon(n) \mbox{ ~ (experiment)} \label{e.1} \\
&& x_k(m+2) = x_k(m) + g_k(n) \gamma(n) \mbox{ ~ (update)}  \label{e.2}
\end{eqnarray}
where
\begin{eqnarray} 
&& g_k(n) = \frac{f(\bx(m + 1)) - f(\bx(m) + \eta_k(n)}{x_k(m + 1) - x_k(m)};~~~~~~~~~~ \label{e.3} \\
&& a_k(n) \mbox{ are independent with} \nonumber \\
&&~~~ P(a_k(n) = -1) = P(a_k(n) = 1) = 0.5 ;\label{e.4}\\
&& \eta_k(n) \mbox{ are independent, zero-mean, bounded}; \label{e.5} \\
&& \epsilon(n), \gamma(n) \in (0, 1), \sum \frac{\gamma^2(n)}{\epsilon^2(n)} < \infty, \sum \gamma(n)\epsilon^2(n) < \infty, \nonumber \\
&&~~~\mbox{ and } \frac{\gamma(n)}{\epsilon^2(n)} \mbox{ is bounded}.  \label{e.6} \\
&& (\mbox{For instance, } \gamma(n) = n^{-0.75}, \epsilon(n) = n^{-0.2}.) \nonumber  
\end{eqnarray}

Our objective is to prove the following theorem.

\begin{thm}
One has
\[
\bx_n \to \bx^*, \mbox{ almost surely as } n \to \infty
\]
where $\bx^*$ is the maximizer of $f(\cdot)$ on $\Re^K$.
\end{thm}

\section{Proof Outline}

Let $\bz(m)$ be the vector with components $z_k(m) = x_k(m) - a_k(n) \epsilon(n) 1\{m = T_k(n) + 1\}$.
That is, $z_k(m)$ is the latest updated value of $x_k$ by time $m$.  Note that $z_k(m)$ does not change when
user $k$ performs an experiment, only when he updates his variable. Of course $\bx$ changes during experiments
and updates, and the gradient is estimated by observing $f(\bx)$, not $f(\bz)$.

Fix any $\delta > 0$. It is shown in the Appendix  that  $u(n) = \| \bz(n \tau) - \bx^*\|$ satisfies the following
two inequalities: 
\begin{equation} \label{e.12}
u(n+1) \leq u(n) -  [\gamma(n) \beta - \alpha(n)], \mbox{ whenever } u(n) > \delta
\end{equation}
and
\begin{equation} \label{e.13}
u^2(n+1) \leq u^2(n) + c(n), \mbox{ whenever } u(n) \leq \delta.
\end{equation}
In these expressions, $\beta > 0, c(n) \to 0,$ and $\sum \alpha(n)$ converges
to a finite random variable.

The claim is that these inequalities imply that $u(n) \leq 3 \delta$ for all $n \geq n_0$ for some finite $n_0$. To see this,
choose $n_0$ so that $c(n) \leq 3 \delta^2$ for $n \geq n_0 - 1$ and
$\sum_{n = n_0}^{n_0 + m} \alpha(n)  \leq \delta$ for all $m \geq 0$.  Let
$n_1$ be the first time after $n_0$ that $u(n) > \delta$. If there is no such time,
we are done.  Else, let $m_1$ be the first time after $n_1$ that $u(n) \leq \delta$.
Such a time must exist because of (\ref{e.12}), for otherwise $u(n) \to - \infty$ since $\sum \gamma(n) = \infty$
and $\sum \alpha(n) < \infty$.
Let then $n_2$ be the first time after $m_1$ that $u(n) > \delta$, then
$m_2$ the first time after $n_2$ that $u(n) \leq \delta$, and so on.
Finally, let $v(j)$ be the maximum value of $u(n)$ for $n \in \{n_j, \ldots, m_j - 1\}$.
Because of (\ref{e.13}), $u^2(n_j) \leq u^2(n_j - 1) + c(n_j - 1) \leq \delta^2 + 3 \delta^2$, so that $u(n_j) \leq 2 \delta$.
Also, because of (\ref{e.12}),
$v(j) - u(n_j) \leq \max_m \sum_{n = n_0}^{n_0 + m} \alpha(n) \leq \delta$.  Hence $v(j) \leq 3 \delta$
for all $j$, so that $u(n) \leq 3 \delta$ for all $n \geq n_0$. 

Since $\delta > 0$ is arbitrary, it follows that $u(n) \to 0$.
Since $\|\bz(n\tau) - \bx(n\tau)\| \leq \epsilon(n)$, this implies that $\bx(n) \to \bx^*$.

\section{Conclusions}

This paper proves the convergence of a distributed version of the Kiefer-Wolfowitz algorithm under some strong
assumptions.  The function is assumed to be strictly concave in a neighborhood of its maximizer and with bounded derivatives up to the third order.
The observation noise is assumed to be bounded.  The agents update periodically, with the same period, but
possibly with different phases.  The proof is self-contained and does not require assuming that the variables
visit a compact set infinitely often. Instead, it shows that the
updates prevent the variables from drifting away.

Many of these assumptions are stronger than necessary.  For instance, the
periods of the different agents could be different.  This assumption can probably be relaxed further by assuming only that
the rates of update converge.  Convergence in probability should occur if only moments of the noise are bounded.
Relaxing the assumptions and a projection version of the algorithm are left for further study.

\section{Acknowledgements}

This work was motivated by an application to wireless networks studied with Piotr  Gawlowicz and Adam Wolisz.
They identified the importance of asynchronous distributed updates in that application. I am grateful for 
their suggestions for this paper.

\section{Appendix:  Proof of (\ref{e.12})-(\ref{e.13})} 

We first 
give the main steps that lead to the inequalities.
The rest of the appendix provides the details of the calculations.

\subsection*{Main Steps}

Inequality (\ref{e.12}) says that when $\bz$ is away from the maximizer $\bx^*$ the gradient updates
bring it closer. This is intuitive since the gradient is then large. 
Inequality (\ref{e.13}) says that when $\bz$ is close to the maximizer, the updates
do not make it move far aways. This happens because the gradient is then small.

Every $\tau$ steps, each variable $x_k$ gets updated roughly in the direction of the partial derivative of $f(\cdot)$
with respect to that variable. Thus, $\bz$ gets updated roughly in the direction of the gradient $\nabla f(\bx)$.
Errors occur because of corruptions of the gradient estimate due to observation noise and the changes of the
other variables by other agents.
More precisely, using (\ref{e.3}) one finds (see Lemma \ref{L4})
\begin{eqnarray}
&&  \bw(n) := \bz(n\tau + \tau) - \bz(n \tau)  \nonumber \\
&&~= \gamma(n) \nabla f(\bz(n\tau)) + \mu(n) \gamma(n)/\epsilon(n) + O(\rho(n)) \nonumber \\
&&~~  \label{e.200} 
\end{eqnarray}
where $\rho(n) = \max\{ \gamma^2(n)/\epsilon(n), \gamma(n)\epsilon^2(n)\}$ and $\mu(n)$ is a bounded random vector that is zero-mean given ${\cal F}_{n-1}$ where
\[
{\cal F}_n := \{a_k(m), \eta_k(m), m \leq n; k = 1, \ldots, K\}.
\] 
Also, in (\ref{e.200}), $O(\rho(n))$ is a random vector whose components are bounded in absolute value by a constant times $\rho(n)$.

Identity (\ref{e.200}) implies (see Lemma \ref{L5}),
\begin{equation} \label{e.201}
\| \bw(n) \|^2 =  O(\gamma^2(n)/\epsilon^2(n)). 
\end{equation}
Hence,
\begin{eqnarray*}
u^2(n+1) &=& \| \bz(n \tau + \tau)  - \bx^* \|^2 = \| \bz(n \tau) - \bx^* + \bw(n) \|^2 \\
&=& u^2(n) + 2 (\bz(n \tau) - \bx^*)' \bw(n) + \| \bw(n)\|^2 \\
&=& u^2(n) + 2 \gamma(n) (\bz(n \tau) - \bx^*)'  \nabla f(\bz(n\tau)) \\
&&~~~~ + 2 (\gamma(n)/\epsilon(n)) (\bz(n \tau) - \bx^*)' \mu(n) \\
&&~~~~ +   2 (\bz(n \tau) - \bx^*)' O(\rho(n))+  O(\gamma^2(n)/\epsilon^2(n)).
\end{eqnarray*}

Now, 
\begin{equation} \label{e.402}
(\bz(n\tau) - \bx^*)' \nabla(f(\bz(n\tau)) \leq f(\bz(n\tau)) - f(\bx^*),
\end{equation}  by concavity. (See Lemma \ref{L6}.)
Thus,
\begin{eqnarray} 
u^2(n+1) &\leq& u^2(n) + 2 \gamma(n) (f(\bz(n \tau)) - f(\bx^*)) \nonumber \\
&&~+ 2 (\gamma(n)/\epsilon(n)) (\bz(n \tau) - \bx^*)' \mu(n)  \nonumber \\
&& ~ + 2 (\bz(n \tau) - \bx^*)'  O(\rho(n)) \nonumber \\
&&~ +   O(\gamma^2(n)/\epsilon^2(n)). \label{e.50}
\end{eqnarray}

When $u(n) > \delta$, one has 
\begin{equation} \label{e.403}
f(\bz(n \tau)) <  f(\bx^*) - \beta u(n)
\end{equation}
for some $\beta > 0$, by
the strict concavity of $f(\cdot)$ around $\bx^*$.  (See Lemma \ref{L7}.) 
Also, 
\[
(\bz(n \tau) - \bx^*)' \mu(n) = u(n) \sum_k h_k(n) \mu_k(n)
\]
with 
\[
h_k(n) = \frac{z_k(n\tau) - x^*_k}{u(n)}.
\]
Hence, when $u(n) > \delta$,
\begin{eqnarray*}
&& u^2(n+1) \leq u^2(n) - 2 \gamma(n) \beta u(n)  \\
&&~~~ + 2u(n) (\gamma(n)/\epsilon(n)) \sum_k h_k(n) \mu_k(n) \\
&&~~~ + 2 u(n) O(\rho(n)) + O(\gamma^2(n)/\epsilon^2(n))\\
&\leq& u^2(n) - 2 \gamma(n) \beta u(n) \\
&&~~~ + 2u(n)  (\gamma(n)/\epsilon(n)) \sum_k h_k(n) \mu_k(n)  \\
&&~~~ + 2u(n) O(\rho(n))  + 2u(n) O(\gamma^2(n)/\epsilon^2(n))/(2\delta) \\
&=& u^2(n) - 2 \gamma(n) \beta u(n) + 2u(n)  (\gamma(n)/\epsilon(n)) \sum_k h_k(n) \mu_k(n)  \\
&&~~~ + 2u(n) [O(\rho(n)) + O(\gamma^2(n)/\epsilon^2(n))/(2\delta)] \\
&\leq& u^2(n) - 2 \gamma(n) \beta u(n) \\
&&~~~ + 2u(n)[ (\gamma(n)/\epsilon(n)) \sum_k h_k(n) \mu_k(n) + O(\kappa(n))]
\end{eqnarray*}

where 
\begin{eqnarray*}
\kappa(n) &=& \max\{\gamma^2(n)/\epsilon^2(n),\rho(n))\} \\
&=& \max\{\gamma^2(n)/\epsilon^2(n), \gamma(n) \epsilon^2(n)\}.
\end{eqnarray*}

Hence,
\begin{equation} \label{e.800}
u^2(n+1)  \leq u^2(n) - 2u(n)[\beta \gamma(n) - \alpha(n)]
\end{equation}
where
\[
\alpha(n) =  (\gamma(n)/\epsilon(n)) \sum_k h_k(n) \mu_k(n) + O(\kappa(n)).
\]
Now, (\ref{e.800}) implies implies (\ref{e.12}), i.e.,
\[
u(n+1) \leq u(n) - [\beta \gamma(n) - \alpha(n)].
\]
Indeed, if this last inequality were violated, one would have
\begin{eqnarray*}
u^2(n+1) &>& \{u(n) - [\beta \gamma(n) - \alpha(n)]\}^2 \\
&=& u^2(n) - 2u(n)[\beta \gamma(n) - \alpha(n)] + [\beta \gamma(n) - \alpha(n)]^2\\
&\geq& u^2(n) - 2u(n)[\beta \gamma(n) - \alpha(n)].
\end{eqnarray*}
and this would contradict (\ref{e.800}). 

To show that $\sum \alpha(n)$ converges to a finite random variable in Lemma \ref{L8}, one uses
the martingale convergence theorem for the first term and the fact that $\sum_n \kappa(n) < \infty$ by (\ref{e.6}).
For the first term, the key observation  is that $h^2_k(n) \leq 1$.  (See Lemma \ref{L8}.)

When $u(n) \leq \delta$, (\ref{e.50}) 
\begin{eqnarray*}
u^2(n+1) &\leq& u^2(n)  + 2 (\gamma(n)/\epsilon(n)) (\bz(n \tau) - \bx^*)' \mu(n)  + O(\kappa(n))\\
&=&  u^2(n) + c(n)
\end{eqnarray*}
with
\[
 c(n) =  2 (\gamma(n)/\epsilon(n)) (\bz(n \tau) - \bx^*)' \mu(n)  +  O(\kappa(n)).
\]
The martingale convergence theorem implies that the first term goes to zero, because $\|\bz(n \tau) - \bx^*\| \leq \delta$
and $\sum \gamma^2(n)/\epsilon^2(n) < \infty$. The last term also goes to zero.
(See Lemma \ref{L9}.)

\vskip 0.1in
The next section develops some estimates.

\subsection*{Preliminary Calculations}

We recall the following notation that avoids having to keep track of explicit constants.

\begin{definition}
Let $\{h(n), n \geq 0\}$ be a sequence of positive numbers.  By definition, $\{O(h(n)), n \geq 0\}$ designates a sequence of
random variables such that 
\[
|O(h(n))| \leq C h(n), \forall n
\]
for some constant $C$.

The same notation is used when the variables $O(h(n))$ are deterministic and in the vector case when the inequality holds componentwise.

\end{definition}

This definition leads immediately to the following properties. (The last one assumes $A > 0, B > 0$ and uses 
\[
O(\gamma(n)/\epsilon(n)) \leq O(\epsilon(n))
\]
since $\gamma(n)/\epsilon^2(n)$ is bounded, by (\ref{e.6}).)

\begin{lem} \label{L1}
One has
\begin{eqnarray}
&& \left[ O(h(n)) \right]^\alpha  = O(h(n)^\alpha), \forall \alpha > 0 \label{e.14} \\
&& O(h_1(n)) \times O( h_2(n)) = O(h_1(n)h_2(n)) \label{e.15} \\
&& O(h_1(n)) + O(h_2(n)) \nonumber \\
&&~~~ = O(\max\{h_1(n), h_2(n)\}) \label{e.16} \\
&&  \mbox{If }h_1(n) \leq C_1 h_2(n) \leq C_2 h_1(n), n \geq n_0,  \nonumber \\
&&~~~  \mbox{ then } O(h_1(n)) = O(h_2(n)) \label{e.17}  \\
&& \max\{\epsilon(n), \gamma(n)(A + B/\epsilon(n))\} = O(\epsilon(n)).\label{e.18} 
\end{eqnarray}
\end{lem}

\begin{lem} \label{L2}

Let $m = T_k(n)$. We claim that
\begin{eqnarray}
&&f(\bx(m+1)) - f(\bx(m)) = a_k(n)\epsilon(n) f_k(\bz(n\tau)) \nonumber \\
&&~~ + V \epsilon(n) + a_k(n) U \epsilon^2(n) + a_k(n) V' \gamma(n)  \nonumber \\
&&~~  +  O(\epsilon^3(n)) \label{e.22c}
\end{eqnarray}
where $U, V, V'$ are bounded and independent of $a_k(n)$ and ${\cal F}_{n-1}$, and $U$ is zero-mean.
Also, $f_k(\bz(n\tau))$
is the partial derivative of $f(\cdot)$ with respect to $x_k$ evaluated at $\bz(n\tau)$.
\end{lem}

\proof{Proof of Lemma \ref{L2}}

Let $m = T_k(n)$. Recall that $T_l(n) + 1$ is the experiment time of agent $l$ during $\{n\tau, \ldots, n\tau + \tau - 1\}$,
so that $T_l(n) + 2$ is his update time. The update equations (\ref{e.1}) and (\ref{e.2}) imply that, 
for $m \in \{n\tau, \ldots, n\tau + \tau - 1\}$,
\[
q_l := x_l(m) - z_l(n\tau) = \left\{
\begin{array}{l l}
a_l(n)\epsilon(n), & \mbox{ if } T_l(n) = T_k(n)  - 1\\
g_l(n)\gamma(n),& \mbox{ if } T_l(n) \leq T_k(n) - 2\\
0, & \mbox{ otherwise}
\end{array}
\right.
\]
and
\[
r_l := x_l(m+1) - z_l(n\tau) = \left\{
\begin{array}{l l}
a_l(n)\epsilon(n), & \mbox{ if } T_l(n) = T_k(n) \\
g_l(n)\gamma(n),& \mbox{ if } T_l(n) \leq T_k(n) - 1\\
0, & \mbox{ otherwise}.
\end{array}
\right.
\]
An important observation is that the gradient estimates $g_l(n)$ for $l \neq k$ are only
affected by $a_k(n)$ at and after time $m + 1$ and then used to update $x_l$ at or after time $m+2$.
Thus, the random variables $q_l(n)$ and $r_l(n)$ for $l \neq k$ that enter in the calculations of $\bx(m)$
and $\bx(m+1)$ are independent of $a_k(n)$.  Moreover,
$a_k(n)$ is independent of $\bz(n\tau)$.

Definition (\ref{e.3}) shows that, for all $l = 1, \ldots, K$,
\[
|g_l(n)| \leq L + G/\epsilon(n) = O(1/\epsilon(n))
\]
where $L$ is the Lipschitz constant and $G$ is the bound on $\eta_k(n)$.

The identities
above show that $\|\br\| = O(\epsilon(n))$ and $\|\bq\| = O(\epsilon(n))$ because $g_l(n) \gamma(n) = O(\epsilon(n))$ by (\ref{e.18}).  
Taylor's theorem implies the following identity:
\begin{eqnarray*}
&& f(\bx(m+1) )- f(\bz(n\tau)) = f(\bz(n\tau) + \br) - f(\bz(n\tau)) \\
&&~~~~ = \br' \nabla f(\bz(n\tau) )  + \frac{1}{2} \br' H \br + O(\epsilon^3(n))
\end{eqnarray*}
where $H = Hf(\bz(n\tau) )$ is the Hessian of $f(\cdot)$ evaluated at $\bz(n\tau)$.

Similarly,
\begin{eqnarray*}
&& f(\bx(m) )- f(\bz(n\tau)) = f(\bz(n\tau) + \bq) - f(\bz(n\tau)) \\
&&~~~~ = \bq' \nabla f(\bz(n\tau) )  + \frac{1}{2} \bq' H \bq + O(\epsilon^3(n)).
\end{eqnarray*}

Subtracting these two expressions, we find
\begin{eqnarray*}
&& f(\bx(m+1) ) - f(\bx(m) ) = (\br - \bq)' \nabla f(\bx(n\tau) ) \\
&&~~~~ + \frac{1}{2} ( \br - \bq)' H  ( \br + \bq) + O(\epsilon^3(n)). 
\end{eqnarray*}

Now,
\[
r_l - q_l = \left\{
\begin{array}{l l}
a_l(n)\epsilon(n), & \mbox{ if } T_l(n) = T_k(n)\\
g_l(n)\gamma(n) - a_l(n)\epsilon(n) ,& \mbox{ if } T_l(n) = T_k(n) - 1\\
0, & \mbox{ otherwise}
\end{array}
\right.
\]
and
\[
r_l + q_l =  \left\{
\begin{array}{l l}
a_l(n)\epsilon(n), & \mbox{ if } T_l(n) = T_k(n)\\
g_l(n)\gamma(n) + a_l(n)\epsilon(n), & \mbox{ if } T_l(n) = T_k(n) - 1\\
2 g_l(n)\gamma(n), & \mbox{ if } T_l(n) \leq  T_k(n) - 2\\
0, & \mbox{ otherwise}.
\end{array}
\right.
\]

In the rest of this proof, $U, U_1,U_2, U_3$ designate random variables that are bounded, zero-mean and
independent of $a_k(n)$ and $V, V', V_1,V_2, V_3, V_4$ designate random variables that are bounded and
independent of $a_k(n)$.

By examining the terms in $\br - \bq$, we finds that
\begin{eqnarray*}
&& (\br - \bq)' \nabla f(\bx(n\tau) ) = a_k(n)\epsilon(n) f_k(\bz(n\tau)) + W
\end{eqnarray*}
where $W$ is a sum of terms of the forms
\[
a_l(n)\epsilon(n) f_l(\bz(n\tau)) \mbox{ and } g_l(n)\gamma(n) f_l(\bz(n\tau)). 
\]
Thus, the terms of the above two types are either of the form
\[
U_1 \epsilon(n) \mbox{ or } V_1 \gamma(n)/\epsilon(n).
\]
We conclude that 
\[
(\br - \bq)' \nabla f(\bx(n\tau) ) = a_k(n)\epsilon(n) f_k(\bz(n\tau))  + U_1 \epsilon(n) + V_1 \gamma(n)/\epsilon(n).
\]

The sum $( \br - \bq)' H  ( \br + \bq)$ is composed of terms that are multiples of one of the following three expressions:
\[
a_i(n)a_j(n) H_{i,j} \epsilon^2(n), a_i(n)g_j(n) H_{i,j}\epsilon(n) \gamma(n), g_i(n)g_j(n) H_{i,j}\gamma^2(n).
\]

Terms of first type yield a sum $a_k(n) U_2 \epsilon^2(n) + V_2 \epsilon^2(n)$ where
$a_k(n) U_2 \epsilon^2(n) = 2 \sum_{j \neq k} a_k(n) a_j(n) H_{k, j}$ and 
\[
V_2 \epsilon^2(n) = H_{k,k} \epsilon^2(n) + \sum_{i \neq k}\sum_{j \neq k} a_i(n)a_j(n) \epsilon^2(n) H_{i, j}.
\]

Terms of the second or third type yield a sum  $U_3 \gamma(n) + a_k(n) V_3 \gamma(n) + V_4 \gamma^2(n)/\epsilon^2(n)$.

Combining the observations above, we conclude that
\begin{eqnarray*}
&& f(\bx(m+1) ) - f(\bx(m) ) \\
&&~ = a_k(n)\epsilon(n) f_k(\bz(n\tau)) + U_1 \epsilon(n) + V_1 \gamma(n)/\epsilon(n) \\
&&~~~ + a_k(n) U_2 \epsilon^2(n)  + V_2 \epsilon^2(n)\\
&&~~~ + U_3 \gamma(n) + a_k(n) V_3 \gamma(n) + V_4 \gamma^2(n)/\epsilon^2(n) + O(\epsilon^3(n)) \\
&&~ = a_k(n)\epsilon(n) f_k(\bz(n\tau)) + V \epsilon(n) + a_k(n) U \epsilon^2(n) \\
&&~~~  + a_k(n) V' \gamma(n) +  O(\epsilon^3(n))
\end{eqnarray*}
where $U, V, V'$ are defined as
\begin{eqnarray*}
&&V \epsilon(n)  = U_1 \epsilon(n)  + V_1 \gamma(n)/\epsilon(n) + V_2 \epsilon^2(n) + U_3 \gamma(n) \\
&&~~~~~ + V_4 \gamma^2(n)/\epsilon^2(n) \\
&& a_k(n) U \epsilon^2 (n) = a_k(n) U_2 \epsilon^2(n) \\
&& a_k(n) V' \gamma(n) =  a_k(n) V_3 \gamma(n).
\end{eqnarray*}

This is (\ref{e.22c}).

You will note that in this derivation, all the terms involving $g_l(n)$ are due to the asynchronous updates where some agents
update while others are estimating the partial derivatives.

\begin{lem} \label{L3}

Let $m = T_k(n)$.  We claim that
\begin{equation}  \label{e.22}
g_k(n) = f_k(\bx(n\tau)) + \mu_k(n)/\epsilon(n) + O(\gamma(n)/ \epsilon(n)) + O(\epsilon^2(n))
\end{equation}
where $\mu_k(n)$ is a bounded random variable that is zero-mean given ${\cal F}_{n-1}$.

\end{lem}

\proof{Proof of Lemma \ref{L3}}

Since $a_k(n) = 1/a_k(n)$ (because $a_k(n) \in \{-1, 1\}$) one has, using Lemma \ref{L2},
\begin{eqnarray*}
&& g_k(n) = \frac{f(\bx(m+1)) - f(\bx(m)) + \eta_k(n) }{a_k(n) \epsilon(n)} \\
&&~= (a_k(n)/\epsilon(n)) \times  [ a_k(n)\epsilon(n) f_k(\bz(n\tau)) + V \epsilon(n) \\
&&~~+ a_k(n) U \epsilon^2(n) + a_k(n) V' \gamma(n) +  O(\epsilon^3(n)) + \eta_k(n) ] \\
&&~= f_k(\bz(n\tau)) +  a_k(n) V + U \epsilon(n) + V'\gamma(n)/\epsilon(n) \\
&&~~ + a_k(n) \eta_k(n)/\epsilon(n) + O(\epsilon^2(n)).
\end{eqnarray*}
This expression is of the form (\ref{e.22}), with 
\[
\mu_k(n)/\epsilon(n) = a_k(n) V + U \epsilon(n) + a_k(n) \eta_k(n)/\epsilon(n)
\]
and
\[
O(\gamma(n)/\epsilon(n)) + O(\epsilon^2(n)) = V' \gamma(n)/\epsilon(n)  + a_k(n) O(\epsilon^2(n)).
\]

\endproof

\subsection*{Proofs of the Main Steps}
\vskip 0.1in

The following Lemma shows that (\ref{e.200}) holds.

\begin{lem} \label{L4}
Let $\bw(n) = \bz(n\tau + \tau) - \bz(n\tau)$.  One has
\begin{equation} \label{e.202}
\bw(n) = \gamma(n) \nabla f(\bz(n\tau)) + (\gamma(n)/\epsilon(n) )\mu(n) + O(\rho(n))
\end{equation}
where $\mu(n)$ is a bounded  random vector that is zero-mean  given ${\cal F}_{n-1}$ and
$\rho(n) = \max\{\gamma^2(n)/\epsilon(n), \gamma(n)\epsilon^2(n)\}$.
\end{lem}
\proof{Proof of Lemma \ref{L4}}

One has
\begin{equation} \label{e.700}
z_k(n\tau + \tau) = z_k(n\tau) + \gamma(n) g_k(n),
\end{equation}
so that Lemma \ref{L3} implies that 
\begin{eqnarray*}
\bw(n) &=& \gamma(n) \nabla f(\bz(n\tau)) + (\gamma(n)/\epsilon(n) )\mu(n) \\
&&~~~+ \gamma(n) [O(\gamma(n)/\epsilon(n)) + O(\epsilon^2(n))] \\
&=& \gamma(n) \nabla f(\bz(n\tau)) + (\gamma(n)/\epsilon(n) )\mu(n) + O(\rho(n)).
\end{eqnarray*}
Hence, (\ref{e.202}) holds.

\endproof

The following Lemma shows that (\ref{e.201}) holds.

\begin{lem} \label{L5}
Let $\bw(n) = \bz(n\tau + \tau) - \bz(n\tau)$.  One has
\begin{equation} \label{e.401}
\| \bw(n) \|^2 = O(\gamma^2(n)/\epsilon^2(n)).  
\end{equation}
\end{lem}

\proof{Proof of Lemma \ref{L5}}
In (\ref{e.200}), which is also (\ref{e.202}), the gradient of $f(\cdot)$ is bounded and so is the vector $\mu(n)$.
Hence,
\[
\bw(n) =  O(\gamma(n)) + O(\gamma(n)/\epsilon(n)) + O(\rho(n)) = O(\gamma(n)/\epsilon(n)).
\]
Thus, (\ref{e.401}) holds.
\endproof

The following Lemma proves (\ref{e.402})

\begin{lem} \label{L6}
One has
\begin{equation} \label{e.406}
(\bz(n\tau) - \bx^*)' \nabla(f(\bz(n\tau)) \leq f(\bz(n\tau)) - f(\bx^*).
\end{equation}
\end{lem}

\proof{Proof of Lemma \ref{L6}}
Let $\bz = \bz(n\tau)$. For $\rho \in [0, 1]$, one has
\[
(1 - \rho) f(\bz) + \rho f(\bx^*) \leq f(\rho \bx^* + (1 - \rho) \bz),
\] 
by concavity of $f(\cdot)$.  By Taylor's theorem,
\[
f(\rho \bx^* + (1 - \rho) \bz) = f(\bz) + \rho (\bx^* - \bz)' \rho f(\bz) + O(\rho^2).
\]
Hence,
\[
f(\bz) + \rho (f(\bx^*) - f(\bz))) \leq f(\bz) + \rho (\bx^* - \bz)' \nabla f(\bz) + O(\rho^2),
\]
so that
\[
\rho (f(\bx^*) - f(\bz)) \leq \rho (\bx^* - \bz)' \nabla f(\bz) + O(\rho^2),
\]
Dividing by $\rho$, we get
\[
f(\bx^*) - f(\bz) \leq (\bx^* - \bz)' \nabla f(\bz) + O(\rho).
\]
Letting $\rho \to 0$ yields (\ref{e.406}).
\endproof

The following Lemma proves  (\ref{e.403}).

\begin{lem} \label{L7}
For any $\delta > 0$, there is some $\beta > 0$ such that
\begin{equation} \label{e.403b}
f(\bz) \leq f(\bx^*) - \beta \|\bz - \bx^*\|, \mbox{ if } \|\bz - \bx^* \| \geq \delta.
\end{equation}
\end{lem}

\proof{Proof of Lemma \ref{L7}}
 By continuity and strict concavity in a neighborhood of $\bx^*$,
 \[
- \alpha :=  \max \{ f(\bz) - f(\bx^*)  \mid \| \bz - \bx^* \| \geq \delta \} < 0.
\]
Let $\beta = \alpha / \delta$.  Assume $\| \bz - \bx^* \| \geq \delta$. Define $\bv$ as follows:
\[
\bv = \rho \bz + (1 - \rho ) \bx^*  \mbox{ with } 1 - \rho = \delta/\| \bz - \bx^* \|.
\]
Then, 
\[
\| \bv - \bx^*\| = \| (1 - \rho) \bz - (1 - \rho) \bx^*\| = (1 - \rho) \|\bz - \bx^*\| = \delta.
\]
Consequently, 
\[
f(\bv) - f(\bx^*) \leq - \alpha.
\]
Also, by concavity,
\[
f(\bv) \geq \rho f(\bx^*) + (1 - \rho) f(\bz).
\]
Hence,
\[
\theta f(\bx^*) + (1 - \rho) f(\bz) \leq f(\bx^*) - \alpha,
\]
so that 
\[
f(\bz) \leq f(\bx^*) - \frac{\alpha}{1 - \rho} = f(\bx^*) - \beta \|\bz - \bx^*\|,
\]
as claimed.

\endproof

The following lemma shows that the sequence $a(n)$ in (\ref{e.12}) sums to a finite random variable.
\begin{lem} \label{L8}
Let
\[
\alpha(n) =  (\gamma(n)/\epsilon(n)) \sum_k h_k(n) \mu_k(n) + O(\kappa(n))
\]
where $\kappa(n) = \max\{\gamma^2(n)/\epsilon^2(n), \gamma(n) \epsilon^2(n)\}$.

Then the sum of $\alpha(n)$ converges to a finite random variable.
\end{lem}

\proof{Proof of Lemma \ref{L8}}
First consider 
\[
 (\gamma(n)/\epsilon(n)) h_k(n) \eta_k(n). 
 \]
 Recall that $|h_k(n)| \leq 1$ and that the random variables $\eta_k(n)$ are bounded and zero-mean given
 ${\cal F}_{n-1}$.  Thus, the sum
 \[
 \sum_{n = 0}^m  (\gamma(n)/\epsilon(n)) h_k(n) \mu_k(n)
 \]
 is a martingale with respect to that filtration ${\cal F}_m$.  Moreover,
 \[
 \sum_n (\gamma(n)/\epsilon(n))^2 < \infty
 \]
by assumption.  Consequently, by the martingale convergence theorem, this sum converges to a finite random variable.
 
 Also, the terms $O( \kappa(n)))$ sum to finite numbers, by (\ref{e.6}).
 \endproof
 
 The following lemma shows that the $c(n)$ in (\ref{e.13}) converge to zero.
 
 \begin{lem} \label{L9}
 Let 
 \[
 c(n) =  2 (\gamma(n)/\epsilon(n)) (\bz(n \tau) - \bx^*)' \mu(n) +  O(\kappa(n))
 \]
 for $n$ such that $\|\bz(n \tau) - \bx^*\| \leq \delta$ and $c(n) = 0$ otherwise.  Then $c(n) \to 0$.
 \end{lem}
 
\proof{Proof of Lemma \ref{L9}}
 Consider the term
 \[
  (\gamma(n)/\epsilon(n))  (z_k(n \tau) - x^*_k) \mu_k(n).
  \]
  Note that
  \begin{eqnarray*}
&&   |  (\gamma(n)/\epsilon(n))  (z_k(n \tau) - x^*_k) |^2 \leq (\gamma^2(n)/\epsilon^2(n)) |z_k(n \tau) - x^*_k|^2 \\
&&~~~~~~~ \leq (\gamma^2(n)/\epsilon^2(n)) \|\bz(n\tau) - \bx^*\|^2 \leq  (\gamma^2(n)/\epsilon^2(n))\delta^2.
\end{eqnarray*}
Consequently, as in Lemma \ref{L8}, these terms sum to a finite random variable.  Hence, the terms
converge to zero.

The terms $O(\kappa(n))$ also converge to zero.

\endproof


\bibliographystyle{nonumber}

 \end{document}